\pdfoutput=1
\RequirePackage{ifpdf}
\ifpdf 
\documentclass[pdftex]{sigma}
\else
\documentclass{sigma}
\fi

\begin{document}

\allowdisplaybreaks

\renewcommand{\thefootnote}{$\star$}

\renewcommand{\PaperNumber}{027}

\FirstPageHeading

\ShortArticleName{The Structure of Line Bundles over Quantum Teardrops}

\ArticleName{The Structure of Line Bundles\\
over Quantum Teardrops\footnote{This paper is a~contribution to the Special Issue on Noncommutative Geometry and Quantum
Groups in honor of Marc A.~Rief\/fel.
The full collection is available at
\href{http://www.emis.de/journals/SIGMA/Rieffel.html}{http://www.emis.de/journals/SIGMA/Rieffel.html}}}

\Author{Albert Jeu-Liang SHEU}

\AuthorNameForHeading{A.J.L.~Sheu}

\Address{Department of Mathematics, University of Kansas, Lawrence, KS 66045, USA}
\Email{\href{mailto:asheu@ku.edu}{asheu@ku.edu}}

\ArticleDates{Received October 07, 2013, in f\/inal form March 15, 2014; Published online March 22, 2014}

\Abstract{Over the quantum weighted 1-dimensional complex projective spaces, called quantum teardrops, the quan\-tum
line bundles associated with the quantum principal ${\rm U}(1)$-bundles introduced and studied by Brzezinski and Fairfax
are explicitly identif\/ied among the f\/initely generated projective modules which are classif\/ied up to isomorphism.
The quantum lens space in which these quantum line bundles are embedded is realized as a~concrete groupoid
$C^*$-algebra.}


\Classification{46L85; 58B32}

\begin{flushright}
\textit{Dedicated to Prof.~Marc A.~Rieffel on the occasion of his 75th birthday}
\end{flushright}

\renewcommand{\thefootnote}{\arabic{footnote}}
\setcounter{footnote}{0}

\section{Introduction}

In the theory of noncommutative topology or geometry~\cite{Conn}, a~generally noncommutative $C^*$-algebra $\mathcal{A}$ or
a~dense ``core'' $*$-subalgebra $\mathcal{A}^{\infty}$ of it is viewed respectively as
the algebra $C(X_{q})$ of continuous functions or the algebra $\mathcal{O}(X_{q})$ of
coordinate functions on an imaginary spatial object $X_{q}$, called a~noncommutative space or a~quantum space.
In many interesting cases, such an imaginary nonexistent space $X_{q}$ is closely related to or actually originated from
a~classical counterpart, a~well-def\/ined topological space or manifold $X$, and we view $X_{q}$ or its
``function algebra'' $C(X_{q})$ or $\mathcal{O}(X_{q})$ as a~quantization of
the classical spatial object $X$.

There have been very intriguing discoveries that a~lot of topological or geometric concepts or properties of a~space $X$
are also carried by (the function algebra of) its quantum counterpart $X_{q}$.
For example, the concept of a~vector bundle $E$~\cite{Huse} over a~compact space $X$ can be reformulated in the
noncommutative context as a~f\/initely generated projective left modules $\Gamma(E_{q})$ over $C(X_{q})$,
viewed as the space of continuous cross-sections of some imaginary noncommutative or quantum vector
bundle $E_{q}$ over $X_{q}$, as suggested by Swan's work~\cite{Swan}.
Beyond the well-known $K$-theoretic study of such noncommutative vector bundles up to stable isomorphism, the
classif\/ication of them up to isomorphism for $C^*$-algebras was made popular by Rief\/fel~\cite{Ri:dsr,Ri:ct} and completed
for some interesting quantum spaces by him and others~\cite{Bach,Pete, Ri:ct,Ri:pm,Sh:ct}.

When the spatial object $X$ is actually a~topological group $G$, the quantization encompasses the group structure by
requiring $C(G_{q})$ or $\mathcal{O}(G_{q})$ to have an additional Hopf $*$-algebra structure,
and we call $G_{q}$ or its function algebra a~quantum group.
Generalizing further, we view a~surjective Hopf $*$-algebra homomorphism $\mathcal{O}(G_{q})\rightarrow\mathcal{O}(H_{q})$
as giving a~quantum subgroup $H_{q}$ of a~quantum group $G_{q}$, and view
the coinvariant $*$-subalgebra $\mathcal{O}(G_{q} /H_{q})$ of $\mathcal{O}(G_{q})$ for the
ca\-no\-nical coaction $\mathcal{O}(G_{q})\overset{\Delta_{R}}{\rightarrow}\mathcal{O}(G_{q})\otimes\mathcal{O}(H_{q})$
as def\/ining a~``quantum homogeneous space''~$G_{q}/H_{q}$.
More generally, given a~coaction $\Delta_{R}:\mathcal{O} (X_{q})\rightarrow\mathcal{O}(X_{q})\otimes\mathcal{O}(H_{q})$
of a~compact quantum group~$H_{q}$ on a~compact quantum space~$X_{q}$, the
coinvariant $*$-subalgebra $\mathcal{O}(X_{q}/H_{q})$ of $\mathcal{O}(X_{q})$ def\/ines
a~``quantum quotient space''~$X_{q}/H_{q}$.

Classically some internal structure of a~vector bundle $E$ over a~space $X$ is often carried by a~principal $G$-bundle~$P$ over~$X$ for some structure group $G$ represented on some vector space~$V$ such that $E=P\times_{G}V$.
The concept of quantum principal bundles has evolved and become well developed through years of
study~\cite{BrzeMaji,Haja:sc}.
In a~recent work of Brzezi\'{n}ski and Fairfax~\cite{BrzeFair}, the quantization of weighted 1-dimensional complex
projective spaces ${\rm WP}(k,l)$, called teardrops by Thurston, and of principal bundles over them was studied.
In particular, the quantum principal ${\rm U}(1)$-bundles and the associated quantum line bundles over the
quantum teardrops ${\rm WP}_{q}(k,l)$ were introduced and analyzed by Brzezi\'{n}ski and Fairfax.
More concretely, they found a~family of quantum line bundles $\mathcal{L}[n] $, $n\in\mathbb{Z}$, inside
a~quantum principal ${\rm U}(1)$-bundle $C(L_{q}(l;1,l))$ over ${\rm WP}_{q}(k,l)$
and showed that the continuous function $C^*$-algebra $C({\rm WP}_{q}(k,l))$ is isomorphic to the
unitization $(\mathcal{K}^{l})^{+}$ of $l$ copies of the algebra $\mathcal{K}$ of compact operators.

In this paper, we f\/irst show that each of $C(L_{q}(l;1,l))$ and $C({\rm WP}_{q}(1,l))$
can be realized as a~concrete groupoid $C^*$-algebra~\cite{Rena}, following the groupoid approach to
study $C^*$-algebras as initiated by Renault~\cite{Rena} and popularized by Curto, Muhly, and Renault~\cite{CuMu,MuRe}.
Then we explicitly identify the completed quantum line bundles $\mathcal{L}[n]$ among the well-known
classif\/ied isomorphism classes of all f\/initely generated projective left modules over $(\mathcal{K}^{l})^{+}$.
This identif\/ication exhibits an interesting connection between ``winding numbers'' and
``ranks''.

\section{Projective modules}\label{section1}

From the analysis point of view, since the category of isomorphism classes of unital commutative $C^*$-algebras is
equivalent to the category of homeomorphism classes of compact Hausdorf\/f spaces, the category of isomorphism classes of
$C^*$-algebras provides a~natural context for the development of noncommutative topology or geometry.

In this context, Swan's theorem~\cite{Swan} makes it legitimate to call an (isomorphism class of) f\/initely generated
projective left module $E$ over a~unital $C^*$-algebra $\mathcal{A}$ an (isomorphism class of) noncommutative vector bundle
over $\mathcal{A}$ or more precisely the (generally imaginary, nonexistent) underlying quantum space.
On the other hand, a~projection $p$ in the $C^*$-algebra $M_{n}(\mathcal{A})$ def\/ines a~left
$\mathcal{A}$-module endomorphism $\xi\in\mathcal{A}^{n}\mapsto\xi p\in\mathcal{A}^{n}$ on the left free
$\mathcal{A}$-module $\mathcal{A}^{n}$, and its image is a~f\/initely generated projective left $\mathcal{A}$-module
$E:=\mathcal{A}^{n}p$.
It is well-known that this association establishes a~bijective correspondence between the unitary equivalence classes of
projections $p$ in $M_{\infty}(\mathcal{A}):=\bigcup_{n=1}^{\infty}M_{n}(\mathcal{A})$, where
$M_{n}(\mathcal{A})$ is embedded in $M_{n+1}(\mathcal{A})$ in the canonical way for each $n$,
and the isomorphism classes of f\/initely generated projective left modules $E$ over $\mathcal{A}$~\cite{Blac}.

Two f\/initely generated projective left modules $E$, $F$ over $\mathcal{A}$ are called stably isomorphic if they become
isomorphic after being augmented by the same f\/initely generated free $\mathcal{A}$-module, i.e.\
$E\oplus \mathcal{A}^{k}\cong F\oplus\mathcal{A}^{k}$ for some $k\in\mathbb{N}$.
The $K_{0}$-group of $\mathcal{A}$ classif\/ies such f\/initely generated projective modules up to stable isomorphism.
The cancellation problem dealing with whether two stably isomorphic f\/initely generated projective left modules are
actually isomorphic goes beyond $K$-theory and is in general an interesting but dif\/f\/icult question.
It was Rief\/fel's pioneering work~\mbox{\cite{Ri:dsr,Ri:ct}} that brought the cancellation problem to the attentions and
interest of researchers in the theory of $C^*$-algebras.
Over some basic geometrically motivated quantum spaces, the f\/initely generated projective left modules have been
successfully classif\/ied~\cite{Bach,Pete, Ri:ct,Ri:pm,Sh:ct}.

As a~simple example, we now describe the classif\/ication of f\/initely generated projective left modules over a~fairly
elementary $C^*$-algebra, which is relevant to our main result later.

Let $\mathcal{K}$ be the algebra of all compact operators on a~separable inf\/inite-dimensional Hilbert space
$\mathcal{H}$, say, $\ell^{2}$.
Recall that for a~$C^*$-algebra $\mathcal{A}$, we use $\mathcal{A}^{+}$ to denote its unitization, a~unital $C^*$-algebra
equal to $\mathcal{A}\oplus\mathbb{C}$ as a~vector space and endowed
with the algebra multiplication $(a,s)(b,t):=(ab+sb+ta,st)$
and involution $(a,s)^{\ast}:=(a^{\ast},\overline{s})$ for $(a,s),(b,t)\in\mathcal{A}\oplus\mathbb{C}$.
In particular, $(\mathcal{K}^{l})^{+}\equiv(\oplus_{s=1}^{l}\mathcal{K})^{+}$ for
$l\in\mathbb{N}$ denotes the unitization of the direct sum of $l$ copies of~$\mathcal{K}$.

The classif\/ication of all isomorphism classes of f\/initely generated projective left modules over~$(\mathcal{K}^{l})^{+}$,
or equivalently, all unitary equivalence classes of projections in $M_{\infty}((\mathcal{K}^{l})^{+})$ is fairly
well understood as summarized below.
In the following, we use~$I$ to denote the multiplicative unit of the unital $C^*$-algebra $(\mathcal{K}^{l})^{+}$,
and $I_{r}$ to denote the identity matrix in $M_{r}((\mathcal{K}^{l})^{+})$, while
\begin{gather*}
P_{n}:=\sum\limits_{i=1}^{n}e_{ii}\in M_{n}(\mathbb{C})\subset \mathcal{K}
\end{gather*}
denotes the standard $n\times n$ identity matrix in $M_{n}(\mathbb{C})\subset\mathcal{K}$ for any integer
$n\geq0$ (with $M_{0}(\mathbb{C})=0$ and $P_{0}=0$ understood).
In particular, $\oplus_{j=1}^{l}P_{k_{j}}\in\mathcal{K}^{l}$ for integers~$k_{j}\geq0$.

\begin{proposition}
The projections $\oplus_{j=1}^{l}P_{k_{j}}\in M_{1}((\mathcal{K}^{l})^{+})$ with $k_{j}
\in\mathbb{Z}_{\geq}:=\big\{ k\in\mathbb{Z}:k\geq0\big\} $ and
\begin{gather*}
I_{r-1}\oplus\big(I-\big({\oplus}_{j=1}^{l}P_{n_{j}}\big) \big) \oplus\big({\oplus}_{j=1}^{l}P_{m_{j}}\big) \in
M_{r+1}\big(\big(\mathcal{K}^{l}\big)^{+}\big)
\end{gather*}
with $r\in\mathbb{N}$ and $n_{j},m_{j}\in\mathbb{Z}_{\geq}$ such that $n_{j}m_{j}=0$ for all $j$ represent all unitarily
inequivalent classes of projections in $M_{\infty}((\mathcal{K}^{l})^{+})$.
\end{proposition}

\section{Quantum spaces and principal bundles}

We recall the def\/inition of a~compact quantum group by Woronowicz~\cite{Woro} as a~unital sepa\-rable $C^*$-algebra
$\mathcal{A}$ with a~comultiplication $\Delta$ such that $\left(\mathcal{A}\otimes1\right)\Delta\mathcal{A}$ and
$\left(1\otimes\mathcal{A}\right)\Delta\mathcal{A}$ are dense in~\mbox{$\mathcal{A}\otimes\mathcal{A}$}.
It is known~\cite{Wo:cm, Woro} that a~compact quantum group $\mathcal{A}$ contains a~dense $*$-sub\-algebra~$\mathcal{A}
^{\infty}$, forming a~Hopf $*$-algebra $\left(\mathcal{A}^{\infty},\Delta,^{\ast},S,\varepsilon\right)$, and has a~Haar
state $h\in \mathcal{A}^{\ast}$ satisfying $h(1)=1$ and
\begin{gather*}
\left(\operatorname{id}\otimes h\right)\Delta a=h(a)1=\left(h\otimes\operatorname{id}\right)\Delta a.
\end{gather*}
We denote $\mathcal{A}^{\infty}$ by $\mathcal{O}\left(G_{q}\right)$ if $\mathcal{A}$ is denoted as $C\left(G_{q}\right)$.

For a~quantum subgroup $H_{q}$ of a~compact quantum group $G_{q}$ given by a~surjective Hopf $*$-algebra homomorphism
$r:\mathcal{O}\left(G_{q}\right)\rightarrow\mathcal{O}\left(H_{q}\right)$, there is a~canonical coaction
$\mathcal{O}\left(G_{q}\right)\overset{\Delta_{R}}{\rightarrow} \mathcal{O}\left(G_{q}\right)\otimes\mathcal{O}\left(H_{q}\right)$
given by $\Delta_{R}:=\left(\operatorname{id}\otimes r\right)\Delta$ for the
comultiplication $\Delta$ of $\mathcal{O}\left(G_{q}\right)$, and the coinvariant $*$-subalgebra
\begin{gather*}
\mathcal{O}\left(G_{q}/H_{q}\right):=\big\{x\in\mathcal{O}\left(G_{q}\right):\Delta_{R}(x)=x\otimes1\big\}
\end{gather*}
for the coaction $\Delta_{R}$ def\/ines a~``quantum homogeneous space''~$G_{q}/H_{q}$.
A fundamental example is the quantum odd-dimensional sphere $S_{q}^{2n+1}={\rm SU}_{q}(n+1)/{\rm SU}_{q}(n)$~\cite{VaSo}
with $q\in(0,1)$ ge\-nerated by $z_{0},\dots ,z_{n}$ subject to the relations
$\sum\limits_{m=0}^{n}z_{m}z_{m}^{\ast} =1$, $z_{i}z_{j}=qz_{j}z_{i}$ for $i<j$, $z_{i}z_{j}^{\ast}=qz_{j}^{\ast} z_{i}$
for $i\neq j$, and $z_{i}z_{i}^{\ast}=z_{i}^{\ast}z_{i}+\left(q^{-2}-1\right)\sum\limits_{m=i+1}^{n}z_{m}z_{m}^{\ast}$.

More generally, given a~coaction $\Delta_{R}:\mathcal{O}(X_{q})\rightarrow\mathcal{O}(X_{q})\otimes\mathcal{O}\left(H_{q}\right)$
of a~compact quantum group~$H_{q}$ on a~compact quantum space $X_{q}$, the
coinvariant $*$-subalgebra
\begin{gather*}
\mathcal{O}\left(X_{q}/H_{q}\right):=\big\{x\in\mathcal{O}(X_{q}):\Delta_{R}(x)=x\otimes1\big\}
\end{gather*}
def\/ines a~``quantum quotient space''~$X_{q}/H_{q}$.
An interesting example is the quantum weighted complex projective space ${\rm WP}_{q}(l_{0},\dots ,l_{n})$ with
$q\in(0,1)$~\cite{BrzeFair}, for pairwise coprime numbers $l_{0},\dots ,l_{n}\in\mathbb{N}$, which is the
quantum quotient space for the coaction
of $\mathcal{O}({\rm U}_{q}(1))\equiv\mathcal{O} ({\rm U}(1))=\mathbb{C}\left[u,u^{\ast}\right]$
on $\mathcal{O}(S_{q}^{2n+1})$ def\/ined by
\begin{gather*}
z_{i}\in\mathcal{O}\big(S_{q}^{2n+1}\big)\mapsto z_{i}\otimes u^{l_{i}}\in\mathcal{O}\big(S_{q}^{2n+1}\big)\otimes\mathcal{O}({\rm U}(1))
\qquad
\text{for}
\quad
i=0,\dots ,n.
\end{gather*}
As special cases, this includes the quantum complex projective space $\mathbb{CP}_{q}^{n}$ when $l_{0}=\dots =l_{n}=1$, and
the so-called quantum teardrop ${\rm WP}_{q}(k,l)$ with coprime $k,l$ when $n=1$.

Brzezi\'{n}ski and Fairfax~\cite{BrzeFair} determined that $S_{q}^{3}$ is a~quantum principal ${\rm U}(1)$-bundle
over ${\rm WP}_{q}(k,l)$, or more precisely, the algebra $\mathcal{O}(S_{q}^{3})$ is a~principal
$\mathcal{O}({\rm U}(1) )$-comodule algebra over $\mathcal{O}({\rm WP}_{q}(k,l) )$, if and only if $k=l=1$.
This result is consistent with the classical ${\rm U}(1)$-action $(z,w)\mapsto(u^{k}z,u^{l}w)$ for $u\in\mathbb{T}$ on $\mathbb{S}^{3}$.
Furthermore they found that the quantum lens space $L_{q}(l;1,l)$~\cite{HongSzym} provides the total space
of a~quantum principal ${\rm U}(1)$-bundle over ${\rm WP}_{q}(1,l)$, where $L_{q}(l;1,l)$ is
the quantum quotient space def\/ined by the coaction $\rho:\mathcal{O}(S_{q}^{3}) \rightarrow\mathcal{O}(S_{q}^{3}) \otimes\mathcal{O}(\mathbb{Z}_{l})$
\qquad
with $\rho(\alpha) =\alpha\otimes w\;$and $\rho(\beta) =\beta\otimes1$ where $\alpha:=z_{0}$ and
$\beta:=z_{1}^{\ast}$ generate $\mathcal{O}(S_{q}^{3}) \equiv\mathcal{O}({\rm SU}_{q}(2))$, and $w$ is the unitary group-like generator
of $\mathcal{O}(\mathbb{Z}_{l})$ with $w^{l}=1$.
More explicitly, $\mathcal{O}(L_{q}(l;1,l) )$ is the $*$-subalgebra of $\mathcal{O}({\rm SU}_{q}(2))$
generated by $c:=\alpha^{l}$ and $d:=\beta$, and a~well-def\/ined coaction
\begin{gather*}
\rho_{l}: \ \mathcal{O}(L_{q}(l;1,l))
\rightarrow\mathcal{O}(L_{q}(l;1,l))\otimes\mathcal{O}({\rm U}(1))
\end{gather*}
with $\rho_{l}(c):=c\otimes u$ and $\rho_{l}(d):=d\otimes u^{\ast}$ makes $\mathcal{O}(L_{q}(l;1,l))$
a~quantum principal ${\rm U}(1)$-bundle over ${\rm WP}_{q}(1,l)$.

Corresponding to the irreducible (1-dimensional) representations of ${\rm U}(1)$ indexed by $n\in\mathbb{Z}$, we
have the irreducible corepresentations of $\mathcal{O}({\rm U}(1))$ on some left comodules denoted
as $W_{n}$.
Following the general theory of constructing f\/initely generated projective modules from quantum principal bundles and
f\/inite-dimensional corepresentations~\cite{BrzeHaja}, Brzezi\'{n}ski and Fairfax took the cotensor product of
$\mathcal{O}(L_{q}(l;1,l))$ with $W_{n}$ over $\mathcal{O}({\rm U}(1))$ to
get a~f\/initely generated projective module $\mathcal{L}[n] \subset\mathcal{O}(L_{q}(l;1,l))$
over $\mathcal{O}({\rm WP}_{q}(1,l))$, naturally called a~quantum line bundle over
${\rm WP}_{q}(1,l)$, and they computed an idempotent matrix $E[n]$ over $\mathcal{O}({\rm WP}_{q}(1,l))$
implementing the projective module $\mathcal{L}[n]$ with complicated
entries $E[n]_{ij}=\omega(u^{n})^{\left[2\right]_{i}}\omega(u^{n})^{\left[1\right]_{j}}$,
where $\omega(u^{n})=\sum\limits_{i}\omega(u^{n})^{\left[1\right]
_{i}}\otimes\omega(u^{n})^{\left[2\right]_{i}}$ comes from a~strong connection
\begin{gather*}
\omega:\ \mathcal{O}({\rm U}(1))\rightarrow \mathcal{O}(L_{q}(l;1,l))
\otimes\mathcal{O} (L_{q}(l;1,l)),
\end{gather*}
and showed in particular that the $\mathcal{O}({\rm WP}_{q}(1,l))$-module $\mathcal{L}[1]$ is not free.

Furthermore Brzezi\'{n}ski and Fairfax found the enveloping $C^*$-algebra of $\mathcal{O}({\rm WP}_{q}(k,l))$
as $C({\rm WP}_{q}(k,l))\cong(\mathcal{K}^{l})^{+}$ and computed its
$K$-groups from the exact sequence
\begin{gather*}
0\rightarrow\mathcal{K}^{l}\equiv\oplus_{j=1}^{l}\mathcal{K}
\rightarrow\big(\mathcal{K}^{l}\big)^{+}\equiv C\left({\rm WP}_{q}\left(k,l\right)\right)\rightarrow\mathbb{C}\rightarrow0.
\end{gather*}

It is then a~natural and interesting question to identify explicitly the completed quantum line bundles
\begin{gather*}
\overline{\mathcal{L}[n]}\equiv C\left({\rm WP}_{q}(1,l)\right)
\otimes_{\mathcal{O}\left({\rm WP}_{q}(1,l)\right)}\mathcal{L}[n]
=\big(\mathcal{K}^{l}\big)^{+}\otimes_{\mathcal{O}\left({\rm WP}_{q}(1,l)\right)}\mathcal{L}[n]
\end{gather*}
over $C({\rm WP}_{q}(1,l))$ for all $n\in\mathbb{Z}$ among the f\/initely generated
projective modules over $(\mathcal{K}^{l})^{+}$ already well classif\/ied.

\section[Quantum lens space as groupoid $C^*$-algebra]{Quantum lens space as groupoid $\boldsymbol{C^*}$-algebra}

In the past, there have been successful studies of the structure of some interesting
$C^*$-algebras~\cite{CuMu,MuRe,SSU,Sh:cqg,Sh:sqs} by realizing them f\/irst as a~concrete groupoid $C^*$-algebra, following the
groupoid approach to $C^*$-algebras initiated by Renault~\cite{Rena} and popularized by the work of Curto, Muhly, and
Renault~\cite{CuMu,MuRe}.
In this section, we f\/irst identify the $C^*$-algebra $C(L_{q}(l;1,l))$ for $q\in(0,1)$
with a~concrete groupoid $C^*$-algebra, and then f\/ind an explicit description of the structure of $C(L_{q}(l;1,l))$.
We construct the groupoid directly from the irreducible representations of $C(L_{q}(l;1,l))$
classif\/ied by Brzezi\'{n}ski and Fairfax~\cite{BrzeFair}.
Our approach should be compared with the machinery developed by Kumjian, Pask, Raeburn, Renault, and Paterson
in~\cite{KPaRaRe,Pa} that associates groupoid $C^*$-algebras to graph $C^*$-algebras.

By Proposition 2.4 of~\cite{BrzeFair}, the faithful $*$-representation $\pi^{\oplus}\equiv\oplus_{s=1}^{l}\pi_{s}$ of
$\mathcal{O}({\rm WP}_{q}(1,l))$ on $\oplus_{s=1}^{l}V_{s}$ factors through the key $*$-representation
$\pi$ of $\mathcal{O}({\rm SU}_{q}(2))$ on $V\cong\oplus_{s=1}^{l}V_{s}$, where each
$V_{s}\cong\ell^{2}(\mathbb{Z}_{\geq})$ for $\mathbb{Z}_{\geq}:=\left\{ k\in\mathbb{Z}:k\geq0\right\} $,
and by Proposition 5.1 of~\cite{BrzeFair}, $\pi^{\oplus}\equiv\oplus_{s=1}^{l}\pi_{s}$ extends to a~faithful
$*$-representation of $C({\rm WP}_{q}(1,l))$ identifying $C({\rm WP}_{q}(1,l))$
with $(\mathcal{K}^{l})^{+}$.

Using the classif\/ication~\cite{BrzeFair} of irreducible $*$-representations
of $\mathcal{O}(L_{q}(l;1,l))\subset\mathcal{O} ({\rm SU}_{q}(2))$ as~$\pi_{s}^{\lambda}$ for $s=0,1,\dots,l$ and
$\lambda\in\mathbb{T}$, we can realize $C(L_{q}(l;1,l))$ as a~groupoid $C^*$-algebra as follows.

For $s>0$ and $\lambda\in\mathbb{T}$, each $\pi_{s}^{\lambda}$ is an irreducible representation of $\mathcal{O}(L_{q}(l;1,l))$
on $\ell^{2}(\mathbb{Z}_{\geq})$ such that $\pi_{s}^{\lambda}(c)$ for any f\/ixed~$s$
is the same weighted unilateral shift independent of~$\lambda$, with strictly positive
weights $\prod\limits_{m=1}^{l}\sqrt{1-q^{2(pl+s-m)}}$ and dif\/ferent from the (backward) unilateral shift~$\mathcal{S}$ on~$\ell^{2}(\mathbb{Z}_{\geq})$, that sends the standard basis vector $e_{p}$ of $\ell
^{2}(\mathbb{Z}_{\geq})$ to $e_{p-1}$ (with $e_{-1}:=0$), only by a~compact operator, while
$\oplus_{s=1}^{l}\pi_{s}^{\lambda}(d)=\lambda(\oplus_{s=1}^{l}\pi_{s}^{1}(d))$
with $\oplus_{s=1}^{l}\pi_{s}^{1}(d)$ a~compact diagonal operator on $\oplus_{s=1}^{l}\ell^{2}(\mathbb{Z}_{\geq})$
with distinct nonzero eigenvalues $q^{pl+s}$, $p\in\mathbb{Z}_{\geq}$.
Applying functional calculus to $\oplus_{s=1}^{l}\pi_{s}^{\lambda}(d)$ to get scaled diagonal matrix units
and then composing with powers of $\oplus_{s=1}^{l}\pi_{s}^{\lambda}(c)$ or its adjoint, we can get all
matrix units for each component $\ell^{2}(\mathbb{Z}_{\geq})$ of $\oplus_{s=1}^{l}\ell^{2}(\mathbb{Z}
_{\geq})$ and hence for each $\lambda\in\mathbb{T}$,
\begin{gather*}
\oplus_{s=1}^{l}\mathcal{K}\left(\ell^{2}(\mathbb{Z}_{\geq})\right)
\subset\overline{\left(\oplus_{s=1}^{l}\pi_{s}^{\lambda}\right)\left(\mathcal{O}\left(L_{q}\left(l;1,l\right)\right)\right)}.
\end{gather*}
On the other hand, $(\oplus_{s=1}^{l}\pi_{s}^{\lambda})(c)$ modulo
$\oplus_{s=1}^{l}\mathcal{K}(\ell^{2}(\mathbb{Z}_{\geq}))$ is the direct sum of $l$ copies of the same unilateral shift~$\mathcal{S}$.
So the image $C^*$-algebra $\overline {\pi_{s}^{\lambda}(\mathcal{O}(L_{q}(l;1,l)))}$
is the standard Toeplitz $C^*$-algebra~$\mathcal{T}$, with $\sigma(\pi_{s}^{\lambda}(c))=\operatorname{id}_{\mathbb{T}}$
and $\sigma(\pi_{s}^{\lambda}(d))=0$ for all~$s$, where
$\sigma:\mathcal{T}\rightarrow C(\mathbb{T} )$ is the standard symbol map of $\mathcal{T}$, while for the
image $C^*$-algebra of the direct sum $\oplus_{s=1}^{l}\pi_{s}^{\lambda}$, we have a~short exact sequence
\begin{gather*}
0\rightarrow\oplus_{s=1}^{l}\mathcal{K}\left(\ell^{2}(\mathbb{Z}_{\geq})\right)
\rightarrow\overline{\left(\oplus_{s=1}^{l}\pi_{s}^{\lambda}\right)\left(\mathcal{O}\left(L_{q}\left(l;1,l\right)\right)\right)}
\rightarrow C\left(\mathbb{T}\right)\rightarrow0.
\end{gather*}

Note that each of the one-dimensional irreducible representations $\pi_{0}^{\mu}$ of $\mathcal{O}(L_{q}(l;1,l))$
with $\pi_{0}^{\mu}(c)=\mu\in\mathbb{T}$ and $\pi_{0}^{\mu}(d)=0$
factors through each $\pi_{s}^{\lambda}$, or more explicitly,
$\pi_{0}^{\mu}=\eta_{\mu}\circ\sigma\circ\pi_{s}^{\lambda}$ for the evaluation character $\eta_{\mu}:C(\mathbb{T})\rightarrow\mathbb{C}$
with $\eta_{\mu}(f):=f(\mu)$.
Hence the $\mathbb{T} $-parameter family $\left\{ \oplus_{s=1}^{l}\pi_{s}^{\lambda}\right\}_{\lambda\in\mathbb{T}}$ of
representations together represent faithfully the enveloping $C^*$-algebra $C(L_{q}(l;1,l))$ of
$\mathcal{O}(L_{q}(l;1,l))$.

More ef\/fectively, we can merge the $\mathbb{T}$-parameter family $\left\{ \oplus_{s=1}^{l}\pi_{s}^{\lambda}\right\}
_{\lambda\in\mathbb{T}}$ of representations of $\mathcal{O}(L_{q}(l;1,l))$ into one
representation $\oplus_{s=1}^{l}\tilde{\pi}_{s}$ on the Hilbert space
$L^{2}(\mathbb{T})\otimes(\oplus_{s=1}^{l}\ell^{2}(\mathbb{Z}_{\geq}))$
or equivalently on $\ell^{2}(\mathbb{Z})\otimes(\oplus_{s=1}^{l}\ell^{2}(\mathbb{Z}_{\geq}))$
via the Fourier transform on $\mathbb{T}$.
More precisely, we have
$\oplus_{s=1}^{l}\tilde{\pi}_{s}(c)=\operatorname{id}_{\ell^{2}(\mathbb{Z})}\otimes
(\oplus_{s=1}^{l}\pi_{s}^{1}(c))$
and $\tilde{\pi}_{s}(d)=\mathcal{U}\otimes(\oplus_{s=1}^{l}\pi_{s}^{1}(d))$ for the (backward) bilateral
shift $\mathcal{U}$ on $\ell^{2}(\mathbb{Z})$.
Clearly $\oplus_{s=1}^{l}\tilde{\pi}_{s}$ is a~faithful representation of $\mathcal{O}(L_{q}(l;1,l))$
and extends to a~faithful representation of $C(L_{q}(l;1,l))$.
In the following, we denote by $\tilde{\pi}^{\oplus}:=\oplus_{s=1}^{l}\tilde{\pi}_{s}$ this faithful representation of
$C(L_{q}(l;1,l))$ on $\ell^{2}(\mathbb{Z})\otimes(\oplus_{s=1}^{l}\ell
^{2}(\mathbb{Z}_{\geq}))$.

Now we consider the ($r$-discrete) groupoid
\begin{gather*}
\mathfrak{G}:=\mathbb{Z}\times\left(\left.
\left(\mathbb{Z}\ltimes\left(\bigsqcup\limits_{s=1}^{l}\mathbb{Z}\right)^{+}\right)
\right\vert_{\left(\bigsqcup\limits_{s=1}^{l}\mathbb{Z}_{\geq}\right)^{+}}\right)
\end{gather*}
which is the direct product of the group $\mathbb{Z}$ and the transformation
groupoid $\mathbb{Z}\ltimes\left(\bigsqcup\limits_{s=1}^{l}\mathbb{Z}\right)^{+}$
restricted to the positive half $\left(\bigsqcup\limits_{s=1}^{l}\mathbb{Z}_{\geq}\right)^{+}$,
where $\left(\bigsqcup\limits_{s=1}^{l}\mathbb{Z}\right)^{+}$
is the one-point compactif\/ication of the disjoint union $\bigsqcup\limits_{s=1}^{l}\mathbb{Z}$ of $l$ copies of $\mathbb{Z}$, and $\mathbb{Z}$
acts canonically
by~translation on each component $\mathbb{Z}$ of $\bigsqcup\limits_{s=1}^{l}\mathbb{Z}\subset\left(\bigsqcup\limits_{s=1}^{l}\mathbb{Z}_{\geq}\right)^{+}$
while f\/ixing the point at inf\/inity $\infty\in\left(\bigsqcup\limits_{s=1}^{l}\mathbb{Z}\right)^{+}$.
More explicitly,
\begin{gather*}
\left(k,m,p\right)_{s}\left(k^{\prime},m^{\prime},p^{\prime}\right)_{s}=\left(k+k^{\prime},m+m^{\prime},p^{\prime}\right)_{s}
\end{gather*}
exactly when $p=p^{\prime}+m^{\prime}$ for $k,k^{\prime},m,m^{\prime} \in\mathbb{Z}$ and
$p,p^{\prime}\in\mathbb{Z}_{\geq}$, where the subscript $s$ in $\left(k,m,p\right)_{s}$ and $\left(k^{\prime},m^{\prime},p^{\prime}\right)_{s}$
indicates that $p$ and $p^{\prime}$ come from the same $s$-th component
of $\bigsqcup\limits_{s=1}^{l}\mathbb{Z}_{\geq}$.
We remark that with the group $\mathbb{Z}^{2}$ being amenable, the full groupoid $C^*$-algebra of $\mathfrak{G}$ is the
same as its reduced groupoid $C^*$-algebra by Proposition 2.15 of~\cite{MuRe}.

Before proceeding further, we introduce an open subgroupoid $\mathfrak{F}$ of $\mathfrak{G}$ def\/ined by
\begin{gather*}
\mathfrak{F}:=\left(\mathbb{Z}\times\left(\left.
\left(\mathbb{Z} \ltimes\left(\bigsqcup\limits_{s=1}^{l}\mathbb{Z}\right)^{+}\right)
\right\vert_{\bigsqcup\limits_{s=1}^{l}\mathbb{Z}_{\geq}}
\right)\right)
\cup\left(\{0\} \times\left(\mathbb{Z}\ltimes\left\{ \infty\right\} \right)\right)\subset\mathfrak{G}.
\end{gather*}

Let $\tilde{\rho}$ be the representation of the groupoid $C^*$-algebra $C^{\ast}(\mathfrak{F})$ induced of\/f
the counting measure $\mu$ supported on the set $\bigsqcup\limits_{s=1}^{l}\{0\} $ that generates the dense
invariant open subset $\bigsqcup\limits_{s=1}^{l}\mathbb{Z}_{\geq}$ of the unit space $\left(\bigsqcup\limits_{s=1}^{l}\mathbb{Z}_{\geq}\right)^{+}$.
By Proposition 2.17 of~\cite{MuRe} (or by a~direct inspection for this fairly simple $r$-discrete groupoid),
$\tilde{\rho}$ is faithful.
We note that the representation space of $\tilde{\rho}$ is isomorphic to
\begin{gather*}
\ell^{2}\left(\mathbb{Z}\times\left(\bigsqcup\limits_{s=1}^{l}\mathbb{Z}_{\geq}\right)\right)
\equiv\ell^{2}\left(\mathbb{Z}\right)\otimes\left(\oplus_{s=1}^{l}\ell^{2}(\mathbb{Z}_{\geq})\right),
\end{gather*}
and that
\begin{gather*}
\tilde{\pi}^{\oplus}\left(c\right)
=\tilde{\rho}\left(\sum\limits_{s=1}^{l}\left(\sum\limits_{p=1}^{\infty}
\left(\prod\limits_{m=1}^{l}\sqrt{1-q^{2\left(pl+s-m\right)}}\right)\delta_{\left(0,-1,p\right)_{s}}\right)\right),
\end{gather*}
where the argument of $\tilde{\rho}$ is understood as an element of $C_{c}(\mathfrak{F})\subset C_{c}\left(\mathfrak{G}\right)$
with value equal to
\begin{gather*}
\lim_{p\rightarrow\infty}\left(\prod\limits_{m=1}^{l}\sqrt{1-q^{2\left(pl+s-m\right)}}\right)=1
\qquad
\text{(for any $s$)}
\end{gather*}
at the point $\left(0,-1,\infty\right)\in\mathfrak{F}\subset\mathfrak{G}$ while vanishing at $\left(k,m,\infty\right)\in\mathfrak{G}$
for all $\left(k,m\right)\neq\left(0,-1\right)$.
Also we have
\begin{gather*}
\tilde{\pi}^{\oplus}\left(d\right)
=\tilde{\rho}\left(\sum\limits_{s=1}^{l}\left(\sum\limits_{p=0}^{\infty}q^{pl+s}\delta_{\left(-1,0,p\right)_{s}}\right)\right),
\end{gather*}
where the argument of $\tilde{\rho}$ is an element of $C_{c}(\mathfrak{F})\subset C_{c}(\mathfrak{G})$
with value equal to $\lim\limits_{p\rightarrow\infty}q^{pl+s}=0$ (for any $s$) at the point $(k,m,\infty)$ for all $(k,m)$.

Now via $\tilde{\rho}^{-1}\circ\tilde{\pi}$, we can view $c$, $d$ as elements of $C_{c}(\mathfrak{F})\subset
C^{\ast}(\mathfrak{F} )$ and hence view $C(L_{q}(l;1,l))$ as embedded in
$C^{\ast}(\mathfrak{F})$.
Applying functional calculus to $d^{\ast}d$, we can get $\mathbb{C}\delta_{(0,0,p)_{s}}\subset C(L_{q}(l;1,l))$
for all $p\in\mathbb{Z}_{\geq}$ and $1\leq s\leq l$, and then by composing with
$c^{\ast}$ and $d^{\ast}$, we get $\mathbb{C}\delta_{(0,1,p)_{s}}$ and $\mathbb{C}\delta_{(1,0,p)_{s}}$
contained in $C(L_{q}(l;1,l))$ for any $p\in\mathbb{Z}_{\geq}$ and $1\leq
s\leq l$, which generate the convolution $*$-subalgebra
\begin{gather*}
C_{c}\left(\mathfrak{F}\Big|_{\bigsqcup\limits_{s=1}^{l}\mathbb{Z}_{\geq}}\right)
\subset C^{\ast}(\mathfrak{F})\overset{\tilde{\rho}}{\subset}\mathcal{B}
\left(\ell^{2}\left(\mathbb{Z}\right)\otimes\left(\oplus_{s=1}^{l}\ell^{2}(\mathbb{Z}_{\geq})\right)\right).
\end{gather*}
On the other hand, for any $n\in\mathbb{Z}$, the $\left\vert n\right\vert $-th power of $c$ or $c^{\ast}$ provides an
element of $C_{c}\left(\mathfrak{F} \right)$ having a~nonvanishing positive value at every point in
\begin{gather*}
\left\{ \left(0,n,p\right)_{s}:p\in\mathbb{Z}_{\geq},1\leq s\leq l\right\} \cup\left\{ \left(0,n,\infty\right)\right\}
\end{gather*}
while vanishing at all other points of $\mathfrak{F}$.
So the $C^*$-subalgebra $C\left(L_{q}\left(l;1,l\right)\right)$ of $C^{\ast}(\mathfrak{F})$ contains all
elements of $C_{c}\left(\mathfrak{F} \right)$ and hence equals $C^{\ast}(\mathfrak{F})$.

We summarize:

\begin{theorem}
$C\left(L_{q}\left(l;1,l\right)\right)\cong C^{\ast}(\mathfrak{F})$, where $\mathfrak{F}$ is the
topological groupoid
\begin{gather*}
\mathfrak{F}=\left[\mathbb{Z}\times\left(\left.
\left(\mathbb{Z} \ltimes\left(\bigsqcup\limits_{s=1}^{l}\mathbb{Z}\right)^{+}\right)
\right\vert_{\left(\bigsqcup\limits_{s=1}^{l}\mathbb{Z}_{\geq}\right)^{+}}\right)
\right]
\backslash\left[\left(\mathbb{Z}\backslash \{0\} \right)\times\left(\mathbb{Z}\ltimes\left\{ \infty\right\} \right)\right].
\end{gather*}
\end{theorem}

In the general theory of groupoid $C^*$-algebras~\cite{Rena}, open invariant subsets and their complements in the unit
space of a~groupoid give rise respectively to closed ideals and quotients of its groupoid $C^*$-algebra, and under suitable
conditions the association is bijective which broadens a~result of Gootman and Rosenberg~\cite{GoRo} for transformation
groups.

Decomposing the base space $\left(\bigsqcup\limits_{s=1}^{l}\mathbb{Z}_{\geq}\right)^{+}$ of $\mathfrak{F}$ into the open
invariant subspace $\bigsqcup\limits_{s=1}^{l}\mathbb{Z}_{\geq}$ and its closed invariant complement $\left\{ \infty\right\} $,
we get the closed ideal
\begin{gather*}
C^{\ast}\left(\mathfrak{F}\Big|_{\bigsqcup\limits_{s=1}^{l}\mathbb{Z}_{\geq}}\right)
=C^{\ast}\left(\mathbb{Z}\times\left(\bigsqcup\limits_{s=1}^{l}\left(\mathbb{Z}\ltimes\mathbb{Z}\right)|_{\mathbb{Z}_{\geq}}\right)\right)
\cong C\left(\mathbb{T}\right)\otimes\mathcal{K}^{l}
\end{gather*}
of $C^{\ast}(\mathfrak{F})$ and the quotient
\begin{gather*}
C^{\ast}(\mathfrak{F})/C^{\ast}\left(\mathfrak{F}\Big|_{\bigsqcup\limits_{s=1}^{l}\mathbb{Z}_{\geq}}\right)\cong
C^{\ast}\left(\mathbb{Z}\ltimes\left\{\infty\right\}\right)\cong C\left(\mathbb{T}\right),
\end{gather*}
which can be summarized as follows.

\begin{corollary}
There is a~short exact sequence of $C^*$-algebras
\begin{gather*}
0\rightarrow C(\mathbb{T})\otimes\mathcal{K}^{l}\rightarrow C(L_{q}(l;1,l))\rightarrow C(\mathbb{T})
\rightarrow0.
\end{gather*}
\end{corollary}

In fact, from the above analysis, we actually have the following explicit description
\begin{gather*}
C(L_{q}(l;1,l))=\big\{(a_{1},\dots,a_{l})\in\oplus_{s=1}^{l}C(\mathbb{T},\mathcal{T} ):
\;
\sigma\circ a_{1}=\dots=\sigma\circ a_{l}\;\text{constant on}\;\mathbb{T}\big\}
\end{gather*}
in terms of the standard Toeplitz $C^*$-algebra $\mathcal{T}$ and its symbol map $\sigma:\mathcal{T}\rightarrow C(\mathbb{T})$.

\section{Line bundles over quantum teardrops}

In this section, we identify concretely the quantum line bundles $\overline{\mathcal{L}[n]}$
over $C({\rm WP}_{q}(1,l))\cong(\mathcal{K}^{l})^{+}$ for $q\in(0,1)$.
First we recall that the coaction $\rho_{l}$ of $\mathcal{O}({\rm U}(1))$ on $\mathcal{O}(L_{q}(l;1,l))$
gives a~$\mathbb{Z}$-grading of $\mathcal{O}(L_{q}(l;1,l))$
with $c$ of degree $1$ and $d$ of degree $-1$, such that $\mathcal{O}({\rm WP}_{q}(1,l))$ generated
by $b:=cd$ and $a:=dd^{\ast}$ is the degree-$0$ component of $\mathcal{O}(L_{q}(l;1,l))$,
while $\mathcal{L}[n]$ is the degree-$n$ component of $\mathcal{O}(L_{q}(l;1,l))$
for general $n\in\mathbb{Z}$~\cite{BrzeFair}.

Now we introduce a~compatible $\mathbb{Z}$-grading on the convolution $*$-algebra $C_{c}(\mathfrak{F})$,
based on the groupoid structure.
We def\/ine the homogeneous degree-$n$ component as $C_{c}(\mathfrak{F})_{n}:=C_{c}\big(\mathfrak{F}_{n}\big)$ for the open set
\begin{gather*}
\mathfrak{F}_{n}:=\bigsqcup\limits_{s=1}^{l}\left\{ (k,k-n,p)_{s}:p\in\mathbb{Z}_{\geq},
n-p\leq k\in\mathbb{Z}\right\} \cup\left\{ (0,-n,\infty)\right\} \subset\mathfrak{F}.
\end{gather*}
Note that $\mathfrak{F}=\bigsqcup\limits_{n\in\mathbb{Z}}\mathfrak{F}_{n}$
and $C_{c}(\mathfrak{F})=\oplus_{n\in\mathbb{Z}}C_{c}(\mathfrak{F}_{n})$ becomes a~$\mathbb{Z}$-graded $*$-algebra
with $\deg(\delta_{(k,m,p)_{s}})$ $=k-m$.
Furthermore $c\in C_{c}(\mathfrak{F}_{1})$ and $d\in C_{c}(\mathfrak{F}_{-1})$ for the
generators $c,d\in\mathcal{O}(L_{q}(l;1,l))\subset C_{c}(\mathfrak{F})$ of
$\mathcal{O}(L_{q}(l;1,l))$.
So this groupoid $\mathbb{Z}$-grading on $C_{c}(\mathfrak{F})$ when restricted to the $*$-subalgebra
$\mathcal{O}(L_{q}(l;1,l))\subset C_{c}(\mathfrak{F})$ coincides with the
original $\mathbb{Z}$-grading on $\mathcal{O}(L_{q}(l;1,l))$.
So when viewed as elements of $C_{c}(\mathfrak{F})$, the elements of $\mathcal{L}[n]
\subset\mathcal{O}(L_{q}(l;1,l))$ are homogeneous of degree $n$.
That~is
\begin{gather*}
\mathcal{L}[n] \subset C_{c}(\mathfrak{F})_{n}\equiv C_{c}\big(\mathfrak{F}_{n}\big).
\end{gather*}

Also note that $C_{c}(\mathfrak{F})_{0}=C_{c}(\mathfrak{F}_{0})$ where
$\mathfrak{F}_{0}\subset\mathfrak{F}$ consisting of $(0,0,\infty)$ and elements of the form $(m,m,p)_{s}$
with $p,m+p\in\mathbb{Z}_{\geq}$ is an open subgroupoid of $\mathfrak{F}$.
It is clear that the $*$-algebra $\mathbb{Z}$-grading structure on $C_{c}(\mathfrak{F})$ makes each
$C_{c}(\mathfrak{F})_{n}$ a~left $C_{c}(\mathfrak{F})_{0}$-module.

By the analysis already done on $\mathcal{L}[0] =\mathcal{O} ({\rm WP}_{q}(1,l))$
in~\cite{BrzeFair} or a~direct analysis of the ge\-ne\-rators $a$, $b$
of $\mathcal{O}({\rm WP}_{q}(1,l))\equiv\mathcal{L}[0] \subset C_{c}(\mathfrak{F}_{0})$, we get
\begin{gather*}
C_{c}(\mathfrak{F}_{0})\subset C({\rm WP}_{q}(1,l))
=C^{\ast}(\mathfrak{F}_{0})\subset C^{\ast}(\mathfrak{F})\equiv C(L_{q}(l;1,l)).
\end{gather*}
In particular, $C({\rm WP}_{q}(1,l))$ is realized as the groupoid $C^*$-algebra of the subgroupoid
$\mathfrak{F}_{0}$ of $\mathfrak{F}$.

Let $\overline{\mathcal{L}[n]}$ be the completion of $\mathcal{L}[n]$
in $C^{\ast}(\mathfrak{F})=C(L_{q}(l;1,l))$.
In the following, we show that $\overline{\mathcal{L}[n]}$ is a~f\/initely generated projective left module
over $C({\rm WP}_{q}(1,l))\subset C^{\ast}(\mathfrak{F})$, and hence we can make the
canonical identif\/ication
\begin{gather*}
\overline{\mathcal{L}[n]}\equiv C({\rm WP}_{q}(1,l))
\otimes_{\mathcal{O}({\rm WP}_{q}(1,l))}\mathcal{L}[n].
\end{gather*}

It is easy to see that the $\mathcal{O}({\rm WP}_{q}(1,l))$-module structure on $\mathcal{L}[n]$ by left multiplication
in $C(L_{q}(l;1,l))$ is consistent
with the $C_{c}(\mathfrak{F})_{0}$-module structure on $C_{c}(\mathfrak{F})_{n}$ under the embeddings of
$\mathcal{O}({\rm WP}_{q}(1,l))\equiv\mathcal{L} [0] \subset C_{c}(\mathfrak{F})_{0}$
and $\mathcal{L}[n] \subset C_{c}(\mathfrak{F})_{n}$ into
$C^{\ast}(\mathfrak{F})=C(L_{q}(l;1,l))$.

On the other hand, we have $C_{c}(\mathfrak{F})_{n} \subset\overline{\mathcal{L}[n]}\subset
C(L_{q}(l;1,l))\equiv C^{\ast}(\mathfrak{F})$, using our knowledge of the
$\left\vert n\right\vert $-th power of $c$ or $c^{\ast}$ and that $C_{c}(\mathfrak{F}_{0})\subset\overline
{\mathcal{L}[0]}$.
So
\begin{gather*}
\overline{\mathcal{L}[n]}=\overline{C_{c}(\mathfrak{F} )_{n}}\subset C^{\ast}(\mathfrak{F})\end{gather*}
for each $n$.

Let $X_{m}:=\bigsqcup\limits_{s=1}^{l}\left\{(p+m,p)_{s}:p\geq0\right\}\subset\mathbb{Z}\times\left(\bigsqcup\limits_{s=1}^{l} \mathbb{Z}_{\geq}\right)$
for $m\in\mathbb{Z}$, with
\begin{gather*}
\ell^{2}\left(\mathbb{Z}\times\left(\bigsqcup\limits_{s=1}^{l}\mathbb{Z}_{\geq}\right)\right)=
{\displaystyle\bigoplus\limits_{m\in\mathbb{Z}}}  \ell^{2}(X_{m}).
\end{gather*}

Note that for all $m\in\mathbb{Z}$,
\begin{gather*}
\tilde{\rho}(c)\big(\ell^{2}(X_{m})\big),
\;
\tilde{\rho}(d^{\ast})\big(\ell^{2}(X_{m})\big)\subset\ell^{2}(X_{m+1})
\end{gather*}
while
\begin{gather*}
\tilde{\rho}(b)\big(\ell^{2}(X_{m})\big),
\;
\tilde{\rho}(a)\big(\ell^{2}(X_{m})\big)\subset\ell^{2}(X_{m}).
\end{gather*}
More generally, for all $m\in\mathbb{Z}$,
\begin{gather*}
\tilde{\rho}\big(\overline{\mathcal{L}[n]}\big)\big(\ell^{2}(X_{m})\big)
=\tilde{\rho}\big(\overline {C_{c}(\mathfrak{F})_{n}}\big)\big(\ell^{2}(X_{m})\big)\subset\ell^{2}(X_{m+n}).
\end{gather*}

Identifying $(p+m,p)_{s}\in X_{m}$ with $p$ in the $s$-th copy of $\mathbb{Z}_{\geq}$ in
$\bigsqcup\limits_{s=1}^{l}\mathbb{Z}_{\geq}$, we get a~unitary operator
\begin{gather*}
u_{m}: \ \ell^{2}(X_{m})\rightarrow\ell^{2}\left(\bigsqcup\limits_{s=1}^{l}\mathbb{Z}_{\geq}\right)\cong\oplus_{s=1}^{l}\ell^{2}
(\mathbb{Z}_{\geq})
\end{gather*}
that intertwines $\tilde{\rho}(b)|_{\ell^{2}(X_{m})}$
and $\tilde{\rho}(a)|_{\ell^{2}(X_{m})}$ with $\pi^{\oplus}(b)$ and $\pi^{\oplus}(a)$ respectively.
More generally, the operator
$u_{m} \circ\tilde{\rho}(f)\circ u_{m-n}^{-1}\in\mathcal{B}(\oplus_{s=1}^{l}\ell^{2}(\mathbb{Z}_{\geq}))$
for $f\in\overline{C_{c}(\mathfrak{F})_{n}}\equiv\overline {C_{c}\big(\mathfrak{F}_{n}\big)}$ is independent of $m$, and hence
$\overline{\mathcal{L}[n]}=\overline{C_{c}(\mathfrak{F})_{n}}$ is embedded isometrically into
$\mathcal{B}(\oplus_{s=1}^{l}\ell^{2}(\mathbb{Z}_{\geq}))$ by
$\rho_{n,m}:=u_{m}\circ\tilde{\rho}(\cdot)\circ u_{m-n}^{-1}$ for any $m\in\mathbb{Z}$.
Note that the $\overline{\mathcal{L}[0]}$-module structure on $\overline{\mathcal{L}[n]}$ is
consistent with the $\rho_{0,0}\big(\overline{\mathcal{L}[0]})$-module structure on
$\rho_{n,0}\big(\overline{\mathcal{L} [n]})$ under the embedding $\rho_{n,0}$, where
\begin{gather*}
\rho_{0,0}\big(\overline{\mathcal{L}[0]}\big)
=C({\rm WP}_{q}(1,l))\cong\big({\oplus}_{s=1}^{l} \mathcal{K}\big)^{+}
\equiv\big({\oplus}_{s=1}^{l}\mathcal{K}(\ell^{2}(\mathbb{Z}_{\geq}))\big)^{+}.
\end{gather*}

Furthermore, since $u_{m}\circ\tilde{\rho}(\chi_{C_{n}})\circ u_{m-n}^{-1}=\oplus_{s=1}^{l}\mathcal{S}^{n}$
with $\mathcal{S}$ the backward unilateral shift on $\ell^{2}(\mathbb{Z}_{\geq})$ as def\/ined previously,
for the characteristic function $\chi_{C_{n}}\in C_{c}\big(\mathfrak{F}_{n}\big)$ of the open and compact set
\begin{gather*}
C_{n}:=\left\{ (0,-n,p)_{s}:n\leq p\in\mathbb{Z}_{\geq}\right\}\cup\left\{(0,-n,\infty)\right\} \subset \mathfrak{F}_{n},
\end{gather*}
we have
\begin{gather*}
u_{m}\circ\tilde{\rho}\big(\overline{\mathcal{L}[n]}\big)\circ
u_{m-n}^{-1}=u_{m}\circ\tilde{\rho}\big(\overline{C_{c}\big(\mathfrak{F}_{n}\big)}\big)\circ
u_{m-n}^{-1}=\big({\oplus}_{s=1}^{l}\mathcal{K}\big)+\mathbb{C}\big({\oplus}_{s=1}^{l}\mathcal{S}^{n}\big)
\end{gather*}
which is isomorphic, as a~left $(\oplus_{s=1}^{l}\mathcal{K})^{+}$-module, to
\begin{gather*}
\Big(\big({\oplus}_{s=1}^{l}\mathcal{K}\big)^{+}\oplus\big({\oplus}_{s=1}^{l}\mathcal{K}\big)^{+}\Big)
\Big(I_{1}\oplus\big({\oplus}_{s=1}^{l}P_{n}\big)\Big)
\end{gather*}
if $n\geq0$, and to
\begin{gather*}
\big({\oplus}_{s=1}^{l}\mathcal{K}\big)^{+}\Big(I-\big({\oplus}_{s=1}^{l}P_{-n}\big)\Big)
\end{gather*}
if $n<0$, where we recall that $I_{1}$ denotes the identity matrix in $M_{1}((\oplus_{s=1}^{l}\mathcal{K})^{+})$
while $I$ denotes the identity element of $(\oplus_{s=1}^{l}\mathcal{K} )^{+}$,
and hence $I_{1}\oplus(\oplus_{s=1}^{l}P_{n})\in M_{2}((\oplus_{s=1}^{l}\mathcal{K})^{+})$ while
\begin{gather*}
I-\big({\oplus}_{s=1}^{l}P_{-n}\big)\in\big({\oplus}_{s=1}^{l}\mathcal{K}\big)^{+}=M_{1}\Big(\big({\oplus}_{s=1}^{l} \mathcal{K}\big)^{+}\Big).
\end{gather*}

As summarized below, we have the modules $\overline{\mathcal{L}[n]}$ identif\/ied concretely among the
f\/initely generated projective left modules over $(\mathcal{K}^{l})^{+}$ enumerated earlier in Section~\ref{section1}.

\begin{theorem}
$\overline{\mathcal{L}[n]}$ is isomorphic to the projective left module
over $C({\rm WP}_{q}(1,l))\cong(\mathcal{K}^{l})^{+}$ for $q\in(0,1)$ determined by the projection
$I_{1}\oplus(\oplus_{j=1}^{l}P_{n})\in M_{2}((\mathcal{K}^{l})^{+})$ if
$n\geq0$, and the projection $I-(\oplus_{j=1}^{l}P_{-n})\in M_{1}((\mathcal{K}^{l})^{+})$ if $n<0$.
\end{theorem}

It is interesting to note that this theorem exhibits some kind of an index relation between the ``winding number''~$n$ of the line bundle $\mathcal{L}[n]$ and the ``rank'' of its representative projection $I_{1}\oplus(\oplus_{j=1}^{l}P_{n})$ or $I-(\oplus_{j=1}^{l}P_{-n})$.

Finally, we mention the classif\/ication of isomorphism classes of f\/initely generated projective left modules over the
quantum 3-sphere by Bach~\cite{Bach} which shows that the projections $1\otimes P_{k}$ with $k\geq0$ and $I_{r}$ with
$r\in\mathbb{N}$ represent all unitarily inequivalent classes of projections in $M_{\infty}(C(S_{q}^{3}))$.
In view of this classif\/ication, we observe that
$C(S_{q}^{3})\otimes_{C({\rm WP}_{q}(1,l))}\overline{\mathcal{L}[n]}$
for all $n\in\mathbb{Z}$ is the same rank-1 free module over $C(S_{q}^{3})$,
showing that these non-isomorphic quantum line bundles $\overline{\mathcal{L}[n]}$ over
${\rm WP}_{q}(1,l)$ pull back to the same quantum line bundles over $S_{q}^{3}$ via the quotient map
$S_{q}^{3}\rightarrow {\rm WP}_{q}(1,l)$.
This phenomenon resembles the well-known classical result that the pull-back, to the total space $P$, of a~vector bundle
$P\times_{G}V\rightarrow X$ associated with a~principal $G$-bundle $P\rightarrow X$ for some $G$-vector space $V$ is
always trivial.
In fact, this classical theorem has a~general quantum counterpart~\cite{Haja:priv}.

\subsection*{Acknowledgements}

The author would like to acknowledge the hospitality and support by the National Center for Theoretical Sciences in
Taipei during his 2013 summer visit, and the support by Jonathan Rosenberg's conference grant, NSF-1266158, to attend
the Conference on Noncommutative Geometry and Quantum Groups at the Fields Institute in Toronto, 2013--2014.

\pdfbookmark[1]{References}{ref}
\LastPageEnding


\begin{thebibliography}{99}
\footnotesize\itemsep=0pt

\bibitem{Bach}
Bach K.A., A cancellation problem for quantum spheres, Ph.D. Thesis, University
  of Kansas, Lawrence, 2003.

\bibitem{Blac}
Blackadar B., {$K$}-theory for operator algebras, \textit{Mathematical Sciences
  Research Institute Publications}, Vol.~5, 2nd ed., Cambridge University
  Press, Cambridge, 1998.

\bibitem{BrzeFair}
Brzezi{\'n}ski T., Fairfax S.A., Quantum teardrops, \href{http://dx.doi.org/10.1007/s00220-012-1580-2}{\textit{Comm. Math. Phys.}}
  \textbf{316} (2012), 151--170, \href{http://arxiv.org/abs/1107.1417}{arXiv:1107.1417}.

\bibitem{BrzeHaja}
Brzezi{\'n}ski T., Hajac P.M., The {C}hern--{G}alois character,
  \href{http://dx.doi.org/10.1016/j.crma.2003.11.009}{\textit{C.~R.~Math. Acad. Sci. Paris}} \textbf{338} (2004), 113--116,
  \href{http://arxiv.org/abs/math.KT/0306436}{math.KT/0306436}.

\bibitem{BrzeMaji}
Brzezi{\'n}ski T., Majid S., Quantum group gauge theory on quantum spaces,
  \href{http://dx.doi.org/10.1007/BF02096884}{\textit{Comm. Math. Phys.}} \textbf{157} (1993), 591--638,
  \href{http://arxiv.org/abs/hep-th/9208007}{hep-th/9208007}.

\bibitem{Conn}
Connes A., Noncommutative geometry, Academic Press, Inc., San Diego, CA, 1994.

\bibitem{CuMu}
Curto R.E., Muhly P.S., {$C^*$}-algebras of multiplication operators on
  {B}ergman spaces, \href{http://dx.doi.org/10.1016/0022-1236(85)90062-X}{\textit{J.~Funct. Anal.}} \textbf{64} (1985), 315--329.

\bibitem{GoRo}
Gootman E.C., Rosenberg J., The structure of crossed product {$C^*$}-algebras:
  a proof of the generalized {E}f\/fros--{H}ahn conjecture, \href{http://dx.doi.org/10.1007/BF01389885}{\textit{Invent. Math.}}
  \textbf{52} (1979), 283--298.


\bibitem{Haja:sc}
Hajac P.M., Strong connections on quantum principal bundles, \href{http://dx.doi.org/10.1007/BF02506418}{\textit{Comm.
  Math. Phys.}} \textbf{182} (1996), 579--617, \href{http://arxiv.org/abs/hep-th/9406129}{hep-th/9406129}.


\bibitem{Haja:priv}
Hajac P.M., Private communication.

\bibitem{HongSzym}
Hong J.H., Szyma{\'n}ski W., Quantum lens spaces and graph algebras,
  \href{http://dx.doi.org/10.2140/pjm.2003.211.249}{\textit{Pacific~J. Math.}} \textbf{211} (2003), 249--263.

\bibitem{Huse}
Husemoller D., Fibre bundles, McGraw-Hill Book Co., New York~-- London~--
  Sydney, 1966.

\bibitem{KPaRaRe}
Kumjian A., Pask D., Raeburn I., Renault J., Graphs, groupoids, and
  {C}untz--{K}rieger algebras, \href{http://dx.doi.org/10.1006/jfan.1996.3001}{\textit{J.~Funct. Anal.}} \textbf{144} (1997),
  505--541.

\bibitem{MuRe}
Muhly P.S., Renault J.N., {$C^*$}-algebras of multivariable {W}iener--{H}opf
  operators, \href{http://dx.doi.org/10.2307/1999493}{\textit{Trans. Amer. Math. Soc.}} \textbf{274} (1982), 1--44.

\bibitem{Pa}
Paterson A.L.T., Graph inverse semigroups, groupoids and their
  {$C^*$}-algebras, \textit{J.~Operator Theory} \textbf{48} (2002), 645--662,
  \href{http://arxiv.org/abs/math.OA/0304355}{math.OA/0304355}.

\bibitem{Pete}
Peterka M.A., Finitely-generated projective modules over the
  {$\theta$}-deformed 4-sphere, \href{http://dx.doi.org/10.1007/s00220-013-1724-z}{\textit{Comm. Math. Phys.}} \textbf{321} (2013),
  577--603, \href{http://arxiv.org/abs/1203.6441}{arXiv:1203.6441}.

\bibitem{Rena}
Renault J., A~groupoid approach to {$C^*$}-algebras, \textit{Lecture Notes in
  Mathematics}, Vol.~793, Springer, Berlin, 1980.

\bibitem{Ri:dsr}
Rief\/fel M.A., Dimension and stable rank in the {$K$}-theory of
  {$C^*$}-algebras, \href{http://dx.doi.org/10.1112/plms/s3-46.2.301}{\textit{Proc. London Math. Soc.}} \textbf{46} (1983),
  301--333.

\bibitem{Ri:ct}
Rief\/fel M.A., The cancellation theorem for projective modules over irrational
  rotation {$C^*$}-algebras, \href{http://dx.doi.org/10.1112/plms/s3-47.2.285}{\textit{Proc. London Math. Soc.}} \textbf{47}
  (1983), 285--302.


\bibitem{Ri:pm}
Rief\/fel M.A., Projective modules over higher-dimensional noncommutative tori,
  \href{http://dx.doi.org/10.4153/CJM-1988-012-9}{\textit{Canad.~J. Math.}} \textbf{40} (1988), 257--338.

\bibitem{SSU}
Salinas N., Sheu A.J.L., Upmeier H., Toeplitz operators on pseudoconvex domains
  and foliation {$C^*$}-algebras, \href{http://dx.doi.org/10.2307/1971454}{\textit{Ann. of Math.}} \textbf{130}
  (1989), 531--565.

\bibitem{Sh:ct}
Sheu A.J.L., A cancellation theorem for modules over the group {$C^*$}-algebras
  of certain nilpotent {L}ie groups, \href{http://dx.doi.org/10.4153/CJM-1987-018-7}{\textit{Canad.~J. Math.}} \textbf{39}
  (1987), 365--427.

\bibitem{Sh:cqg}
Sheu A.J.L., Compact quantum groups and groupoid {$C^*$}-algebras,
  \href{http://dx.doi.org/10.1006/jfan.1996.2999}{\textit{J.~Funct. Anal.}} \textbf{144} (1997), 371--393.

\bibitem{Sh:sqs}
Sheu A.J.L., The structure of quantum spheres, \href{http://dx.doi.org/10.1090/S0002-9939-01-06042-7}{\textit{Proc. Amer. Math. Soc.}}
  \textbf{129} (2001), 3307--3311.

\bibitem{Swan}
Swan R.G., Vector bundles and projective modules, \href{http://dx.doi.org/10.1090/S0002-9947-1962-0143225-6}{\textit{Trans. Amer. Math.
  Soc.}} \textbf{105} (1962), 264--277.

\bibitem{VaSo}
Vaksman L.L., Soibelman Ya.S., Algebra of functions on the quantum group {${\rm
  SU}(n+1)$}, and odd-dimensional quantum spheres, \textit{Leningrad Math.~J.}
  \textbf{2} (1991), 1023--1042.

\bibitem{Wo:cm}
Woronowicz S.L., Compact matrix pseudogroups, \href{http://dx.doi.org/10.1007/BF01219077}{\textit{Comm. Math. Phys.}}
  \textbf{111} (1987), 613--665.

\bibitem{Woro}
Woronowicz S.L., Compact quantum groups, in Sym\'etries Quantiques ({L}es
  {H}ouches, 1995), North-Holland, Amsterdam, 1998, 845--884.

\end{thebibliography}
\end{document}